\documentclass[%
  a4paper,
  twocolumn,
]{mypreprint}

\usepackage[english]{babel}

\usepackage{amsmath}
\usepackage{amssymb}
\usepackage{amsthm}

\theoremstyle{definition}\newtheorem{remark}{Remark}

\usepackage{tikz}

\usepackage{pgfplots}
  \pgfplotsset{compat = 1.13}
  \pgfkeys{/pgf/number format/.cd,1000 sep={\,}}
  \usepgfplotslibrary{external}
  \tikzset{external/system call = {%
    pdflatex \tikzexternalcheckshellescape
      -halt-on-error
      -interaction=batchmode
      -jobname "\image" "\texsource"}}
  \tikzexternaldisable

\usepackage{subcaption}

\usepackage[ruled]{algorithm}
\usepackage{algpseudocode}

\usepackage{scalerel}
\renewcommand\widehat[1]{\hstretch{2}{\hat{\hstretch{.5}{#1}}}}

\newcommand{\C}{\ensuremath{\mathbb{C}}}
\newcommand{\R}{\ensuremath{\mathbb{R}}}

\newcommand{\Ar}{\ensuremath{\skew4\widehat{A}}}
\newcommand{\Br}{\ensuremath{\skew4\widehat{B}}}
\newcommand{\Cr}{\ensuremath{\skew4\widehat{C}}}
\newcommand{\Er}{\ensuremath{\skew4\widehat{E}}}
\newcommand{\Gr}{\ensuremath{\skew4\widehat{G}}}
\newcommand{\Hr}{\ensuremath{\skew4\widehat{H}}}
\newcommand{\Kr}{\ensuremath{\skew4\widehat{K}}}
\newcommand{\Mr}{\ensuremath{\skew4\widehat{M}}}

\newcommand{\tx}{\ensuremath{\tilde{x}}}

\DeclareMathOperator{\diag}{diag}

\newcommand{\trans}{\ensuremath{\mkern-1.5mu\mathsf{T}}}

\definecolor{matlabblue}{HTML}{0072BD}
\definecolor{matlaborange}{HTML}{D95319}
\definecolor{matlabpurple}{HTML}{7E2F8E}

\tikzstyle{fom} = [
  black,
  line width = 2.5pt,
  solid
]

\tikzstyle{sovf1} = [
  matlabblue,
  line width = 2pt,
  dashed
]

\tikzstyle{sovf2} = [
  matlaborange,
  line width = 2pt,
  dash dot
]

\algrenewcommand\algorithmicindent{.75\baselineskip}

\begin{document}


\title{Structured vector fitting framework for mechanical systems}

\author[$\ast$]{Steffen W. R. Werner}
\affil[$\ast$]{%
  Courant Institute of Mathematical Sciences, New York University,
  New York, NY 10012, USA.\authorcr
  \email{steffen.werner@nyu.edu}, \orcid{0000-0003-1667-4862}
}

\author[$\dagger$]{Ion Victor Gosea}
\affil[$\dagger$]{%
  Max Planck Institute for Dynamics of Complex Technical Systems,
  Sandtorstr. 1, 39106 Magdeburg, Germany.\authorcr
  \email{gosea@mpi-magdeburg.mpg.de}, \orcid{0000-0003-3580-4116}
}

\author[$\ddagger$]{Serkan Gugercin}
\affil[$\ddagger$]{%
  Department of Mathematics and Division of Computational Modeling and Data
  Analytics, Academy of Data Science, Virginia Tech,
  Blacksburg, VA 24061, USA.\authorcr
  \email{gugercin@vt.edu}, \orcid{0000-0003-4564-5999}
}

\shorttitle{Structured vector fitting framework for mechanical systems}
\shortauthor{S.~W.~R. Werner, I.~V. Gosea, S. Gugercin}
\shortdate{2021-10-17}
\shortinstitute{}

\abstract{%
  In this paper, we develop a structure-preserving formulation of the
  data-driven vector fitting algorithm for the case of modally damped
  mechanical systems.
  Using the structured pole-residue form of the transfer function
  of modally damped second-order systems, we propose two possible structured
  extensions of the barycentric formula of system transfer functions.
  Integrating these new forms within the classical vector fitting algorithm
  leads to the formulation of two new algorithms that allow the computation of
  modally damped mechanical systems from data in a least squares fashion.
  Thus, the learned model is guaranteed to have the desired structure.
  We test the proposed algorithms on two benchmark models.
}


  
\novelty{}

\maketitle


\section{Introduction}%
\label{sec:intro}

Data-driven reduced-order modeling (DD-ROM) is essential in constructing
high-fidelity compact models to approximate the underlying physical phenomena
when an explicit model, a state-space formulation with access to internal
variables, is not available yet an abundant input/output data is.
Thus, DD-ROM circumvents the need to access an exact description of the original
model and is applicable when traditional intrusive projection-based model
reduction is not.
As in the projection case, it is important that the learned model inherits the
physical meaning and structures of the system that has generated the data.    
This is the setup we are interested in here.
Our goal is to develop a data-driven structure-preserving modeling framework
for mechanical systems described by second-order dynamics.

Data in our setting will correspond to transfer function (frequency domain)
samples of the underlying mechanical system.
Let $H(s)$ denote this transfer function and let $\xi_{k}$ denote the sampling
frequencies (points).
Thus, we assume access to the data (measurements) $h_{k} = H(\xi_{k})$, for
$k = 1, 2, \dots, \ell$.
The goal of DD-ROM in this setting is to construct a reduced transfer
function (a low-order rational function) $\Hr(s)$ such that
$\Hr(\xi_{k}) \approx h_{k} = H(\xi_{k})$ in an appropriate measure.
We will call this unstructured (or first-order fitting in this paper) since the
only requirement in this case is that $\Hr(s)$ is a rational function and thus
corresponds to a transfer function with a first-order state-space form.
In this setting, the barycentric rational form of the approximant plays a
crucial role; see~\cite{berrut06barycentric}.
The Loewner framework from~\cite{AA86,mayo2007fsg} that enforces interpolation
of the data, the Vector Fitting (VF) algorithm from~\cite{morGusS99} that
minimizes a least-squares distance, and the AAA algorithm from~\cite{NST18}
that combines interpolation and least squares are just three of the many
techniques for rational data fitting.
We refer the reader to~\cite[Sec.~2.1]{Rod20} for further references. 

Second‐order systems are an important class of structured dynamical systems 
used to describe, for example, the dynamics of mechanical systems, and in
particular, their vibrational response.
Since the underlying second-order structure corresponds to important physical
properties, retaining this structure is vital so that the learned model is
physically meaningful.
Therefor, given the frequency response samples of such systems, our goal is to
construct a \emph{structure-preserving DD-ROM}, in the sense that the learned
model can be interpreted as the transfer function of a second-order (mechanical)
system.
Note that not every rational function can be written as the transfer function
of a second-order system (although the reverse is true).
There have been some recent works on constructing data-driven second-order
models in the interpolatory Loewner framework;
see~\cite{schulze2018data, morPonGB20a}.
There is a much wider literature on projection-based structure-preserving
model reduction for second-order systems.
We refer the reader to~\cite{SaaSW19, Wer21} and the references therein for
details on the projection-based approaches that we do not consider here. 

In this paper, we will focus on structure-preserving second-order DD-ROM using
the least-squares measure.
More specifically, we will enforce modal-damping structure in the learned model.
We will achieve this goal by extending the VF algorithm to the structured
setting.
Up to now, VF has been developed to produce unstructured rational approximants.
We will revise the barycentric formula behind the VF approximant such that
upon convergence the learned model has the desired second-order structure.
This new formulation of the barycentric form will lead to a sequence of linear
least-squares problems whose structure will also inherit the underlying
second-order dynamics.

The rest of the paper is organized as follows:
After providing an overview of the classical VF approach and modally damped
second-order systems in \Cref{sec:background}, we develop the modified
barycentric forms and the resulting structure-preserving VF approaches
together with the corresponding proposed numerical algorithms in
\Cref{sec:sovf}.
The proposed methods are then tested on two benchmark examples in
\Cref{sec:examples}, followed by the conclusions and future research
directions in \Cref{sec:conclusions}.


\section{Background}%
\label{sec:background}

In this section, we provide a brief overview of the classical vector fitting
algorithm and summarize the key structural features of the special class of
mechanical systems under consideration.


\subsection{Classical vector fitting approach}%
\label{subsec:vf}

Assume that one has access to the samples of the transfer function of an
underlying single-input/single-output (SISO) dynamical system to be modeled,
$H(s)$, at the sampling points (frequencies)
$\xi_{1}, \xi_{2}, \ldots, \xi_{\ell} \in \C$.  
Given the data $\{ H(\xi_{i}) \}_{i = 1}^{\ell}$, the goal is to construct
(learn) a degree-$r$ scalar rational function $\Hr(s)$ to solve the
nonlinear rational least-squares (LS) problem
\begin{align} \label{eqn:ls}
  \min\limits_{\Hr} \sum\limits_{i = 1}^{\ell}
    \lvert \Hr(\xi_{k}) - H(\xi_{k}) \rvert^{2}.
\end{align}
Let $\Hr(s) = \frac{n(s)}{d(s)}$ where $n(s)$ and $d(s)$ are, respectively,
degree-$(r-1)$ and degree-$r$ polynomials in $s$.
In other words, $H(s)$ is parametrized by its denominator and numerator
coefficients.
Inserting this form of $\Hr(s)$ into~\cref{eqn:ls}, one can re-write the nonlinear
LS error to minimize as
\begin{align*}
  \sum\limits_{i = 1}^{\ell} \lvert \Hr(\xi_{i}) - H(\xi_{i}) \rvert^{2}
    & = \sum\limits_{i = 1}^{\ell} \frac{1}{\lvert d(\xi_{i}) \rvert^{2}}
    \lvert n(\xi_{i}) - d(\xi_{i})H(\xi_{i})\rvert^{2}.
\end{align*}
The nonlinearity results from the nonlinear dependence of the error on $d(s)$.
To solve this nonlinear LS problem, starting with an initial guess of
$\Hr^{(0)}(s) = \frac{d^{(0)}(s)}{n^{(0)}(s)}$, \cite{SK63} proposed an
iterative scheme where in the $k$-th step the error term~\cref{eqn:ls} is
replaced by
\begin{align} \label{eqn:lsndlin}
  \begin{aligned}
    \sum\limits_{i = 1}^{\ell} \lvert \Hr^{(k)}(\xi_{i}) - H(\xi_{i})
      \rvert^{2} & = \\
    \sum\limits_{i = 1}^{\ell}  \frac{1}{\lvert d^{(k-1)}(\xi_{i}) \rvert^{2}}
      & \lvert n^{(k)}(\xi_{i}) - d^{(k)}(\xi_{i})H(\xi_{i})\rvert^{2}. 
  \end{aligned}
\end{align}
Note that the new error term~\cref{eqn:lsndlin} is now linear in the variables
$n^{(k)}$ and $d^{(k)}(s)$.
Therefor, the SK iteration in~\cite{SK63} converts the original nonlinear LS
problem~\cref{eqn:ls} into solving a sequence of weighted linear LS
problems~\cref{eqn:lsndlin}. 

There are various equivalent forms to represent the rational function $\Hr(s)$.
One can work with numerator and denominator coefficients as the unknowns, or the
poles and residues, for example.
A numerically efficient formulation is the so-called barycentric
representation; see~\cite{berrut06barycentric}.
Let $\Hr^{(k)}(s)$ denote the iterate in the $k$-th step of the SK iteration
as above.
Also let $\lambda_{1}^{(k)}, \ldots, \lambda_{r}^{(k)}$ be mutually distinct
points.
Then, $\Hr^{(k)}(s)$ can be written in the barycentric form as
\begin{align} \label{eqn:barycentric}
  \Hr^{(k)}(s) & = \frac{n^{(k)}(s)}{d^{(k)}(s)}
    = \frac{\sum\limits_{j = 1}^{r} \frac{\phi^{(k)}_{j}}{s -
    \lambda^{(k)}_{j}}}{1 + \sum\limits_{j = 1}^{r}
    \frac{\varphi^{(k)}_{j}}{s - \lambda^{(k)}_{j}}},
\end{align}
where $\{ \phi_{j}^{(k)} \}_{j=1}^{r}$ and $\{ \varphi_{j}^{(k)} \}_{j=1}^{r}$
are the barycentric weights. 
Note that the $\lambda_{j}^{(k)}$'s are not the poles of $\Hr(s)$.
We refer the reader to, e.g.,~\cite{morDrmGB15b} to switch between the
pole-residue form and the barycentric form.

Now inserting $n^{(k)}(s)$ and $d^{(k)}(s)$ from~\cref{eqn:barycentric}
into~\cref{eqn:lsndlin}, in the $k$-th step of the SK iteration, one needs to
solve the weighted linear LS problem 
\begin{align} \label{eqn:vfls}
  \min\limits_{x^{(k)}} \lVert \Delta^{(k)}(A^{(k)} x^{(k)} - h)
    \rVert_{2}^{2},
\end{align}
where
\begingroup
\allowdisplaybreaks
\small
\begin{align} \label{eqn:weightRHS}
  \Delta^{(k)} & = \diag \left( \frac{1}{\lvert d^{(k)}(\xi_{1}) \rvert},
    \ldots, \frac{1}{\lvert d^{(k)}(\xi_{\ell}) \rvert} \right), \quad
    h = \begin{bmatrix} H(\xi_{1}) \\ \vdots \\ H(\xi_{\ell}) \end{bmatrix},\\
  \label{eqn:cauchy}
  A^{(k)} & = \begin{bmatrix}
    \frac{1}{\xi_{1} - \lambda_{1}^{(k)}} & \cdots &
    \frac{1}{\xi_{1} - \lambda_{r}^{(k)}} &
    \frac{-H(\xi_{1})}{\xi_{1} - \lambda_{1}^{(k)}} & \cdots &
    \frac{-H(\xi_{1})}{\xi_{1} - \lambda_{r}^{(k)}} \\
    \vdots & & \vdots & \vdots & & \vdots \\
    \frac{1}{\xi_{\ell} - \lambda_{1}^{(k)}} & \cdots &
    \frac{1}{\xi_{\ell} - \lambda_{r}^{(k)}} &
    \frac{-H(\xi_{\ell})}{\xi_{\ell} - \lambda_{1}^{(k)}} & \cdots &
    \frac{-H(\xi_{\ell})}{\xi_{\ell} - \lambda_{r}^{(k)}}
    \end{bmatrix},
\end{align}
\endgroup
for the solution vector
\begin{align*}
  x^{(k)} & = \begin{bmatrix} \phi_{1}^{(k)} & \cdots & \phi_{r}^{(k)} &
    \varphi_{1}^{(k)} & \cdots & \varphi_{r}^{(k)} \end{bmatrix}^{\trans},
\end{align*}
which forms $\Hr^{(k)}(s)$ at the $k$-th step.
In addition to incorporating the barycentric form into the SK iterations,
\cite{morGusS99} have also observed that the only restrictions on
$\{ \lambda_{j}^{(k)} \}_{j = 1}^{r}$ is to be distinct and they can be updated
at every step.
This is precisely what~\cite{morGusS99} have proposed, leading to the
\emph{Vector Fitting (VF)} algorithm.
VF updates $\{ \lambda_{j}^{(k)} \}_{j = 1}^{r}$ as the roots of denominator
$d^{(k)}(s)$.
Making again use of the barycentric representation~\cref{eqn:barycentric},
these roots are actually the eigenvalues of
$\Ar^{(k)} - \Gr^{(k)} \Cr^{(k)}$, where
\begin{align} \label{eqn:zerosA}
  \Ar^{(k)} & = \diag(\lambda_{1}^{(k)}, \ldots, \lambda_{r}^{(k)}),\\
  \label{eqn:zerosBC}
  \Gr^{(k)} & = \begin{bmatrix} \varphi_{1}^{(k)} ~ \ldots ~
    \varphi_{r}^{(k)} \end{bmatrix}^{\trans} ~~\text{and}~~
  \Cr^{(k)} = \begin{bmatrix} 1 ~ \ldots ~ 1 \end{bmatrix}^{\trans}.
\end{align}
If the algorithm converges, due to the $\{ \lambda_{j}^{(k)} \}_{j = 1}^{r}$
updating strategy, $d^{(k)}(s) \to 1$ and thus the final approximation is
obtained in the pole-residue form with the denominator being $1$
in~\cref{eqn:barycentric}.
The resulting method is summarized in \Cref{alg:vf}, and
we refer the reader to~\cite{morGusS99}, \cite[Chap7]{grivet2015passive}
and~\cite{morDrmGB15} for further details.

\begin{algorithm}[t]
  \caption{(Unstructured) Vector Fitting (VF)}\smallskip
  \label{alg:vf}
  \small
  \hspace*{1.2\baselineskip} \textbf{Input:} Vector $h$ of data
    samples~\cref{eqn:weightRHS}, initial guess for\par
    \hspace{1.2\baselineskip}\phantom{\textbf{Input: }}
    $\{ \lambda_{j}^{(1)} \}_{j = 1}^{r}$. \\
  \hspace*{1.2\baselineskip} \textbf{Output:} Learned ROM
    $\Hr(s) = \Cr (s I_{r} - \Ar)^{-1} \Br$.
  \begin{algorithmic}[1]
    \State Initialize $\Delta^{(1)} = I_{\ell}$ and $k = 1$.
    \While{not converged}
      \State Construct the coefficient matrix $A^{(k)}$ in~\cref{eqn:cauchy}.
      \State Solve the weighted linear least-squares problem~\cref{eqn:vfls}.
      \State Update the expansion points $\{ \lambda_{j}^{(k+1)} \}_{j = 1}^{r}$
        to be the\par
        \hspace{-1.2\baselineskip} eigenvalues of
        $\Ar^{(k)} - \Gr^{(k)} \Cr^{(k)}$
        using~\cref{eqn:zerosA} and~\cref{eqn:zerosBC}.
      \State Update the weighting matrix $\Delta^{(k+1)}$
        by~\cref{eqn:weightRHS}.
      \State Increment $k \leftarrow k + 1$.
    \EndWhile
    \State Set the final ROM matrices to be
      $\Ar = \diag(\lambda_{1}^{(k)}, \ldots, \lambda_{r}^{(k)})$,
      $\Br = \begin{bmatrix} \phi_{1}^{(k)} & \ldots & \phi_{r}^{(k)}
      \end{bmatrix}^{\trans}$ and
      $\Cr = \begin{bmatrix} 1 & \ldots & 1\end{bmatrix}^{\trans}$.
  \end{algorithmic}\smallskip
\end{algorithm}


\subsection{Modally damped second-order systems}%
\label{subsec:mdso}

Next, we take a look at the pole-residue formulation of the structured system
class we consider here, namely the modally damped second-order systems.
As in the previous section, for simplicity we restrict the analysis to the
SISO case.
Assume we have a second-order system of the form
\begin{align*}
  \begin{aligned}
    M \ddot{q}(t) + E \dot{q}(t) + K q(t) & = B_{\mathrm{u}} u(t),\\
    y(t) & = C_{\mathrm{p}} q(t),
  \end{aligned}
\end{align*}
with $M, E, K \in \R^{n \times n}$, $B_{\mathrm{u}} \in \R^{n}$,
$C_{\mathrm{p}}^{\trans} \in \R^{n}$, and modal damping
$E M^{-1} K = K M^{-1} E$ as in~\cite{morBeaB14}.
Note here that for the mechanical system case one additionally
has $M = M^{\trans} > 0$, $E = E^{\trans} \geq 0$, and
$K = K^{\trans} > 0$. This assumption is not necessary in general and
instead we only assume that the pencil $\lambda M - K$ is diagonalizable,
since with  modal damping all three system matrices are simultaneously
diagonalizable.
First, we consider the generalized eigenvalue problems
\begin{align*}
  \begin{aligned}
    K X & = M X \Omega^{2}, &
    K^{\trans} Y & = M^{\trans} Y \Omega^{2},
  \end{aligned}
\end{align*}
where the eigenvector matrices $X$ and $Y$ are scaled such that
\begin{align*}
  \begin{aligned}
    Y^{\trans} M X & = \Omega^{-1}, &
    Y^{\trans} K X & = \Omega,
  \end{aligned}
\end{align*}
with $\Omega = \diag(\omega_{1}, \ldots, \omega_{n})$.
Due to  modal damping, the damping matrix can also be diagonalized
such that
\begin{align*}
  Y^{\trans} E X & = 2 \Psi,
\end{align*}
where $\Psi = \diag(\psi_{1}, \ldots, \psi_{n})$ are the damping ratios of the
system.
Then, the transfer function $H(s)$ satisfies
\begin{align} \nonumber
  H(s) & = C_{\mathrm{p}} (s^{2} M + s E + K)^{-1} B_{\mathrm{u}}\\ \nonumber
  & = C_{\mathrm{p}} X (s^{2} \Omega^{-1} + 2s \Psi + \Omega)^{-1}
    Y^{\trans} B_{\mathrm{u}}\\ \nonumber
  & = \sum\limits_{j = 1}^{n} \frac{\omega_{j} \phi^{\pm}_{j}}
    {s^{2} + 2 \psi_{j} \omega_{j} s + \omega_{j}^{2}}\\ \label{eqn:poleres2}
  & = \sum\limits_{j = 1}^{n} \frac{\omega_{j}\phi^{\pm}_{j}}
    {(s - \lambda_{j}^{+})(s - \lambda_{j}^{-})},
\end{align}
where the pairwise poles of the system are given by 
\begin{align} \label{eqn:quadlambda}
  \lambda^{\pm}_{j} & = -\omega_{j}\psi_{j} \pm \omega_{j}
    \sqrt{\psi_{j}^{2} - 1}.
\end{align}
Since every second-order system can also be written in its first-order
form, we can also write $H(s)$ in the generic  pole-residue formulation as
\begin{align} \label{eqn:poleres1}
  H(s) & = \sum\limits_{j = 1}^{2n}\frac{\phi_{j}}{s - \lambda_{j}} =
    \sum\limits_{j = 1}^{n}\frac{\phi_{j}^{+}}{s - \lambda_{j}^{+}}  + \sum\limits_{j = 1}^{n}\frac{\phi_{j}^{-}}{s - \lambda_{j}^{-}}.
\end{align}
Note that modal damping and the second-order structure enforce additional
properties in the generic pole-residue form, which means that
only for the underlying second-order systems those two formulations,
i.e.,~\cref{eqn:poleres2} and~\cref{eqn:poleres1},
are equivalent. A clear and important advantage of~\cref{eqn:poleres2} is the
enforcement of the underlying system structure.


\begin{table*}
  \vspace{-\baselineskip}
  \small
  \begin{align} \label{eqn:so1A}
    A_{\mathrm{so1}}^{(k)} & = \begin{bmatrix}
      \frac{\omega_{1}^{(k)}}
      {(\xi_{1} - \lambda_{1}^{+,(k)})(\xi_{1} - \lambda_{1}^{-,(k)})}
      & \cdots & 
      \frac{\omega_{r}^{(k)}}
      {(\xi_{1} - \lambda_{r}^{+,(k)})(\xi_{1} - \lambda_{r}^{-,(k)})} &
	  \frac{-H(\xi_{1})}{(\xi_{1} - \lambda_1^{+,(k)})} &
	  \frac{-H(\xi_{1})}{(\xi_{1} - \lambda_1^{-,(k)})}
	  & \cdots &
	  \frac{-H(\xi_{1})}{(\xi_{1} - \lambda_1^{+,(k)})} &
	  \frac{-H(\xi_{1})}{(\xi_{1} - \lambda_{r}^{-,(k)})}
      \\
      \vdots & & \vdots &  \vdots & \vdots & & \vdots & \vdots
      \\
	  \frac{\omega_{1}^{(k)}}
	  {(\xi_{\ell} - \lambda_{1}^{+,(k)})(\xi_{\ell} - \lambda_{1}^{-,(k)})}
	  & \cdots &
	  \frac{\omega_{r}^{(k)}}
	  {(\xi_{\ell} - \lambda_{r}^{+,(k)})(\xi_{\ell} - \lambda_{r}^{-,(k)})} &
	  \frac{-H(\xi_{\ell})}{(\xi_{\ell} - \lambda_{1}^{+,(k)})} &
	  \frac{-H(\xi_{\ell})}{(\xi_{\ell} - \lambda_{1}^{-,(k)})}
	  & \cdots &
	  \frac{-H(\xi_{\ell})}{(\xi_{\ell} - \lambda_{r}^{+,(k)})} &
	  \frac{-H(\xi_{\ell})}{(\xi_{\ell} - \lambda_{r}^{-,(k)})}
	\end{bmatrix}
  \end{align}
  \vspace{.25\baselineskip}
  \hrule
\end{table*}

\section{Second-order vector fitting algorithms}%
\label{sec:sovf}

The classical VF algorithm as outlined in \Cref{subsec:vf} produces an
unstructured rational LS fit.
In this section, we will develop a structured version of VF to model
second-order modally damped system. 
We will achieve this goal by employing the special pole-residue
formulation~\cref{eqn:poleres2} in VF and by modifying the corresponding
barycentric form appearing in VF.
We will propose two formulations for the revised barycentric form and analyze
both forms.
At the end of the newly developed structured VF iteration, the learned model
will be guaranteed to have the modally damped form.


\subsection{Partially structured barycentric form}%
\label{subsec:partstruct}

In our first approach, we develop a second-order VF formulation for modally
damped systems using a partially structured barycentric formulation.
The method computes a second-order system by enforcing $\Hr^{(k)}$, the
re\-duced-order model at iteration step $k$, to have the transfer function 
\begin{align} \label{eqn:partial_bary_Hr}
  \Hr^{(k)}(s) & = \frac{\sum\limits_{j = 1}^{r}
    \frac{\omega_{j}^{(k)} \phi^{\pm,(k)}_{j}}
    {(s - \lambda_{j}^{+,(k)})(s - \lambda_{j}^{-,(k)})}}
    {1 + \sum\limits_{i = 1}^{r} \frac{\varphi_{i}^{+,(k)}}
    {s - \lambda_{i}^{+,(k)}} + \sum\limits_{i = 1}^{r}
    \frac{\varphi_{i}^{-,(k)}}{s - \lambda_{i}^{-,(k)}}},
\end{align}
In other words, the form~\cref{eqn:partial_bary_Hr}
replaces~\cref{eqn:barycentric} in VF. 
The motivation for the revised form~\cref{eqn:partial_bary_Hr} stems from the
desired modally damped structure.
Recall that as classical VF converges, the denominator converges to $1$ and the
numerator becomes the final reduced model.
In the structured form~\cref{eqn:partial_bary_Hr}, we keep the denominator as
before in the classical pole-residue form~\cref{eqn:poleres1}.
However, the numerator is replaced by the structured pole-residue
form~\cref{eqn:poleres2}.
Therefor, upon convergence, the final reduced model, given by the numerator
in~\cref{eqn:partial_bary_Hr}, is guaranteed to have the desired form.

We now discuss the structure of the resulting second-order VF algorithm.
As in \Cref{subsec:vf}, using the relaxation
step~\cref{eqn:lsndlin} we solve a sequence of weighted linear
LS problems of the form
\begin{align} \label{eqn:sovf1ls}
  \min\limits_{\tx^{(k)}} \lVert \Delta^{(k)} (A_{\mathrm{so1}}^{(k)}
    \tx^{(k)} - h) \rVert_{2}^{2},
\end{align}
for the solution vector 
\begingroup
\small
\begin{align*}
  \tx^{(k)} & = \begin{bmatrix} \phi_{1}^{\pm, (k)} ~\cdots~
    \phi_{r}^{\pm, (k)} ~ \varphi_{1}^{+,(k)} ~ \varphi_{1}^{-,(k)} ~ \cdots ~
    \varphi_{r}^{+,(k)} ~ \varphi_{r}^{-,(k)} \end{bmatrix}^{\trans},
\end{align*}
\endgroup
which determines $\Hr^{(k)}(s)$, where $h$ and $\Delta^{(k)}$ are as
in~\cref{eqn:weightRHS}, and the new coefficient matrix
$A_{\mathrm{so1}}^{(k)}$ as in~\cref{eqn:so1A}.
The new coefficient matrix $A_{\mathrm{so1}}^{(k)}$ encodes the underlying
second-order structure. 
Consequently, we replace Steps~3 and~4 in \Cref{alg:vf}
with~\cref{eqn:so1A} and~\cref{eqn:sovf1ls} in the proposed second-order VF
iteration.
Using~\cref{eqn:quadlambda}, the stiffness and damping
coefficients of the pole pairs are given by
\begin{align*}
  \omega_{j}^{(k)} & = \sqrt{\lambda_{j}^{+, (k)} \lambda_{j}^{-, (k)}},\\ 
  \psi_{j}^{(k)} & = -\frac{1}{2 \omega_{j}^{(k)}}
    (\lambda_{j}^{+, (k)} + \lambda_{j}^{-, (k)}),
\end{align*}
for $j = 1, \ldots, r$.
This formulation is needed in constructing $A_{\mathrm{so1}}^{(k)}$
in~\cref{eqn:so1A}, as well as to set up the final data-driven second-order model $\Hr(s) = \Cr_{\mathrm{p}} (s^{2} \Mr + s \Er + \Kr)^{-1}
\Br_{\mathrm{u}}$ where
\begin{align} \label{eqn:learnso}
  \begin{aligned}
    \Mr & = \diag \left( \frac{1}{\omega_{1}^{(k)}}, \ldots, 
      \frac{1}{\omega_{r}^{(k)}} \right),\\
    \Er & = \diag(2 \psi_{1}^{(k)}, \ldots, 2 \psi_{r}^{(k)}),\\
    \Kr & = \diag(\omega_{1}^{(k)}, \ldots, \omega_{r}^{(k)}),\\
    \Br_{\mathrm{u}} & = \begin{bmatrix} \phi_{1}^{\pm, (k)} & \ldots
      & \phi_{r}^{\pm, (k)} \end{bmatrix}^{\trans} \quad\text{and}\\
    \Cr_{\mathrm{p}} & = \begin{bmatrix} 1 & \ldots & 1
      \end{bmatrix}^{\trans}.
  \end{aligned}
\end{align}
A brief sketch of the resulting second-order VF algorithm is given in
\Cref{alg:sovf1}.

\begin{algorithm}[t]
  \caption{Structured Vector Fitting -- Version 1}\smallskip
  \label{alg:sovf1}
  \small
  \hspace*{1.2\baselineskip} \textbf{Input:} Vector $h$ of data
    samples~\cref{eqn:weightRHS}, initial guess for\par
    \hspace{1.2\baselineskip}\phantom{\textbf{Input: }}
    $\{ \lambda_{j}^{+,(1)} \}_{j = 1}^{r}$
    and $\{ \lambda_{j}^{-,(1)} \}_{j = 1}^{r}$. \\
  \hspace*{1.2\baselineskip} \textbf{Output:} Learned ROM
    $\Hr(s) = \Cr_{\mathrm{p}} (s^{2} \Mr + s \Er + \Kr)^{-1}
    \Br_{\mathrm{u}}$.
  \begin{algorithmic}[1]
    \State Initialize $\Delta^{(1)} = I_{\ell}$ and $k = 1$.
    \While{not converged}
      \State Construct the coefficient matrix $A_{\mathrm{so1}}^{(k)}$
        in~\cref{eqn:so1A}.
      \State Solve the weighted linear least-squares
        problem~\cref{eqn:sovf1ls}.
      \State Update the expansion points
        $\{ \lambda_{j}^{\pm,(k+1)} \}_{j = 1}^{r}$ to be the\par
        \hspace{-1.2\baselineskip} eigenvalues of $\Ar^{(k)} - \Gr^{(k)}
        \Cr^{(k)}$ using~\cref{eqn:zerosA} and~\cref{eqn:zerosBC}.
      \State Update the weighting matrix $\Delta^{(k+1)}$
        by~\cref{eqn:weightRHS}.
      \State Increment $k \leftarrow k + 1$.
    \EndWhile
    \State Set the final ROM matrices using~\cref{eqn:learnso}.
  \end{algorithmic}\smallskip
\end{algorithm}

\begin{table*}
  \vspace{-\baselineskip}
  \small
  \begin{align} \label{eqn:so2A}
    A_{\mathrm{so2}}^{(k)} & = \begin{bmatrix}
      \frac{\omega_{1}^{(k)}}
      {(\xi_{1} - \lambda_{1}^{+,(k)})(\xi_{1} - \lambda_{1}^{-,(k)})}
      & \cdots & 
      \frac{\omega_{r}^{(k)}}
      {(\xi_{1} - \lambda_{r}^{+,(k)})(\xi_{1} - \lambda_{r}^{-,(k)})} &
      \frac{-\omega_{1}^{(k)} H(\xi_{1})}
      {(\xi_{1} - \lambda_{1}^{+,(k)})(\xi_{1} - \lambda_{1}^{-,(k)})}
      & \cdots & 
      \frac{-\omega_{r}^{(k)} H(\xi_{1})}
      {(\xi_{1} - \lambda_{r}^{+,(k)})(\xi_{1} - \lambda_{r}^{-,(k)})}
      \\
      \vdots & & \vdots &  \vdots & & \vdots
      \\
      \frac{\omega_{1}^{(k)}}
      {(\xi_{\ell} - \lambda_{1}^{+,(k)})(\xi_{\ell} - \lambda_{1}^{-,(k)})}
      & \cdots & 
      \frac{\omega_{r}^{(k)}}
      {(\xi_{\ell} - \lambda_{r}^{+,(k)})(\xi_{\ell} - \lambda_{r}^{-,(k)})} &
      \frac{-\omega_{1}^{(k)} H(\xi_{\ell})}
      {(\xi_{\ell} - \lambda_{1}^{+,(k)})(\xi_{\ell} - \lambda_{1}^{-,(k)})}
      & \cdots & 
      \frac{-\omega_{r}^{(k)} H(\xi_{1})}
      {(\xi_{\ell} - \lambda_{r}^{+,(k)})(\xi_{\ell} - \lambda_{r}^{-,(k)})}
    \end{bmatrix}
  \end{align}
  \vspace{.25\baselineskip}
  \hrule
\end{table*}

\begin{remark}[Splitting of expansion points]%
  \label{rmk:exppts}
  ~\\
  Another major difference to the classical VF is the splitting of the
  expansion points into two groups $\{ \lambda_{j}^{+, (k)} \}_{j = 1}^{r}$ and
  $\{ \lambda_{j}^{-, (k)} \}_{j = 1}^{r}$, related to each other
  by~\cref{eqn:quadlambda}.
  For mechanical systems with real realizations, the splitting of complex
  points in conjugate pairs with positive imaginary parts
  ($\lambda_{j}^{+, (k)}$) and negative imaginary parts ($\lambda_{j}^{-, (k)}$)
  comes naturally.
  In case of real expansion points, a physics-inspired splitting is with
  respect to bifurcation, i.e., with respect to a centered point on the real
  axis at which the real points would collide and split into complex conjugate
  pairs.
  For simplicity, we assume the real expansion points lie all in the left open
  half-plane.
  Then, we would sort the points such that those with largest magnitude
  ($\lambda_{j}^{-, (k)}$) are paired with those of smallest magnitude
  ($\lambda_{j}^{+, (k)}$).
\end{remark}


\subsection{Fully structured barycentric form}%
\label{subsec:fullstruct}

A second revised barycentric form for $\Hr^{(k)}$ is to replace both the
numerator and denominator by second-order-type pole-residue
forms~\cref{eqn:poleres2}, i.e., we write $\Hr^{(k)}$ as 
\begin{align} \label{eqn:quadbary}
  \Hr^{(k)}(s) & = \frac{\sum\limits_{j = 1}^{r}
    \frac{\omega_{j}^{(k)} \phi^{\pm, (k)}_{j}}
    {(s - \lambda_{j}^{+, (k)})(s - \lambda_{j}^{-, (k)})}}
    {1 + \sum\limits_{j = 1}^{r}
    \frac{\omega_{j}^{(k)} \varphi^{\pm, (k)}_{j}}
    {(s - \lambda_{j}^{+, (k)})(s - \lambda_{j}^{-, (k)})}}.
\end{align}
As in \Cref{subsec:partstruct}, this new barycentric form changes the
form of the weighted linear LS problem in the resulting structured VF
algorithm.
Using~\cref{eqn:quadbary}, we obtain
\begin{align} \label{eqn:sovf2ls}
  \min\limits_{\tx^{(k)}} \lVert \Delta^{(k)} (A_{\mathrm{so2}}^{(k)}
    \tx^{(k)} - h) \rVert_{2}^{2},
\end{align}
where the least-squares matrix $A_{\mathrm{so2}}^{(k)}$ is given
in~\cref{eqn:so2A}, and the weighting matrix and data samples are as
in~\cref{eqn:weightRHS}. 
Then the  solution vector 
\begin{align*}
  \tx^{(k)} & = \begin{bmatrix} \phi_{1}^{\pm, (k)} & \cdots &
    \phi_{r}^{\pm, (k)} & \varphi_{1}^{\pm, (k)} & \cdots &
    \varphi_{r}^{\pm, (k)} \end{bmatrix}^{\trans}
\end{align*}
yields the resulting second-order system  as in~\cref{eqn:learnso}.
The splitting of the expansion points also works as in \Cref{rmk:exppts}.
However, the updating step of the expansion points
(\Cref{alg:sovf1} Step~5) changes.
The denominator of~\cref{eqn:quadbary} corresponds to a second-order system
rather then a first-order system.
While it would be possible to also rewrite this second-order system in
first-order form, the zeros of the denominator are actually given by the
eigenvalues of the quadratic matrix pencil
\begin{align} \label{eqn:quadeig}
  \lambda^{2} \Mr^{(k)} + \lambda \Er^{(k)} +  \left( \Kr^{(k)} +
    \Gr_{\mathrm{u}}^{(k)} \Cr_{\mathrm{p}}^{(k)} \right),
\end{align}
where $\Mr^{(k)}$, $\Er^{(k)}$, $\Kr^{(k)}$ and $\Cr_{\mathrm{p}}^{(k)}$ are
constructed as their final learned counterparts in~\cref{eqn:learnso}, and
\begin{align} \label{eqn:eigdist}
  \Gr^{(k)} & = \begin{bmatrix} \varphi_{1}^{\pm, (k)} & \ldots &
    \varphi_{r}^{\pm, (k)} \end{bmatrix}^{\trans}.
\end{align}
A brief sketch of the resulting method is given in \Cref{alg:sovf2}.

\begin{algorithm}[t]
  \caption{Structured Vector Fitting -- Version 2}\smallskip
  \label{alg:sovf2}
  \small
  \hspace*{1.2\baselineskip} \textbf{Input:} Vector $h$ of data
    samples~\cref{eqn:weightRHS}, initial guess for\par
    \hspace{1.2\baselineskip}\phantom{\textbf{Input: }}
    $\{ \lambda_{j}^{+,(1)} \}_{j = 1}^{r}$
    and $\{ \lambda_{j}^{-,(1)} \}_{j = 1}^{r}$. \\
  \hspace*{1.2\baselineskip} \textbf{Output:} Learned ROM
    $\Hr(s) = \Cr_{\mathrm{p}} (s^{2} \Mr + s \Er + \Kr)^{-1}
    \Br_{\mathrm{u}}$.
  \begin{algorithmic}[1]
    \State Initialize $\Delta^{(1)} = I_{\ell}$ and $k = 1$.
    \While{not converged}
      \State Construct the coefficient matrix $A_{\mathrm{sol2}}^{(k)}$
        in~\cref{eqn:so2A}.
      \State Solve the weighted linear least-squares
        problem~\cref{eqn:sovf2ls}.
      \State Update the expansion points $\{ \lambda_{j}^{\pm, (k+1)}
        \}_{j = 1}^{r}$ to be the\par
        \hspace{-1.2\baselineskip} eigenvalues of~\cref{eqn:quadeig}
        using~\cref{eqn:learnso} and~\cref{eqn:eigdist}.
      \State Update the weighting matrix $\Delta^{(k+1)}$
        by~\cref{eqn:weightRHS}.
      \State Increment $k \leftarrow k + 1$.
    \EndWhile
    \State Set the final ROM matrices using~\cref{eqn:learnso}.
  \end{algorithmic}\smallskip
\end{algorithm}

We note that realness of the resulting state-space realization can be preserved
similar to the classical VF; see~\cite{morGusS99}.
Indeed this task becomes simpler in case of \Cref{alg:sovf2} due to
the natural pairing of complex conjugate expansion points; cf.
\Cref{rmk:exppts}.


\section{Numerical examples}%
\label{sec:examples}

We test the two proposed approaches on two benchmark problems: 
\begin{enumerate}
    \item the butterfly gyroscope example from the~\cite{morwiki_gyro},
    \item the artificial fishtail model from~\cite{SieKM19}.
\end{enumerate}
For simplicity, we consider here only single-input/single-output versions of
these examples.
While the outputs of the first model are summed together, only the second
output entry of the second model is used.
In both examples, we consider data sets with $1000$ linearly equidistant
sampling points on the positive imaginary axis.
For the butterfly gyroscope example, the points lie in
$[10^{2}, 10^{6}]$\,rad/s, while for the artificial fishtail model, the
points are in $[0, 1\,000]$\,rad/s.
Since the matrices of the original mechanical systems are real-valued, the
data samples are closed under conjugation.
This is done by additionally including the complex conjugate counterparts of
both the evaluations and the sampling points into the data sets.
Both proposed structured VF algorithms from \Cref{sec:sovf} are applied
to these two models. 
Thereby, we denote the approach from \Cref{alg:sovf1} using the
partially structured barycentric form by \textsf{SOVF1} and the method in
\Cref{alg:sovf2} based on the fully structured barycentric form by
\textsf{SOVF2}.

The experiments reported here have been executed on a machine equipped with an
AMD Ryzen 5 5500U processor running at 2.10\,GHz and equipped with
16\,GB total main memory.
The computer runs on Windows 10 Home version 20H2 (build 19042.1237) with
MATLAB 9.9.0.1592791 (R2020b).

\begin{center}%
  \setlength{\fboxsep}{5pt}%
  \fbox{%
  \begin{minipage}{.92\linewidth} \small
    \textbf{Code and data availability}\newline
    The source code, authored by Steffen W. R. Werner, of the implementations
    used to compute the presented results, the used data and the computed
    results are available at
    \begin{center}
        \href{https://doi.org/10.5281/zenodo.5539944}%
          {\texttt{doi:10.5281/zenodo.5539944}}
    \end{center}
    under the BSD-2-Clause license.
  \end{minipage}}
\end{center}


\begin{figure}[t]
  \centering
  \begin{subfigure}[b]{\linewidth}
    \centering
  \tikzexternalenable%
  \tikzsetnextfilename{sovf_butterfly_tf}%
  \begin{tikzpicture}
  \pgfplotstableread{graphics/data/sovf1_butterfly.dat}\tableTF
  \pgfplotstableread{graphics/data/sovf2_real_butterfly.dat}\tableTFV
  
  \begin{loglogaxis}[%
    width  = .75\linewidth,
    height = .11\textheight,
    scale only axis,
    xmin = 100,
    xmax = 1000000,
    ymin = 1e-8,
    ymax = 5e-2,
    xminorticks = false,
    yminorticks = false,
    xlabel = {Frequency $\omega$ (rad/s)},
    ylabel = {Magnitude},
    ylabel style   = {yshift = -.3em},
    scaled x ticks = false]
    
    \addplot[fom] table[x index = 0, y index = 1] {\tableTF};
    \addplot[sovf1] table[x index = 0, y index = 2] {\tableTF};
    \addplot[sovf2] table[x index = 0, y index = 2] {\tableTFV};
  \end{loglogaxis}
\end{tikzpicture}%
  \tikzexternaldisable%

    \caption{Transfer functions.}
    \label{fig:sovf_btterfly_tf}
  \end{subfigure}
  \vspace{0\baselineskip}
  
  \begin{subfigure}[b]{\linewidth}
    \centering
  \tikzexternalenable%
  \tikzsetnextfilename{sovf_butterfly_relerr}%
  \begin{tikzpicture}
  \pgfplotstableread{graphics/data/sovf1_butterfly.dat}\tableERR
  \pgfplotstableread{graphics/data/sovf2_real_butterfly.dat}\tableERRV
  
  \begin{semilogyaxis}[%
    xmode  = log,
    width  = .75\linewidth,
    height = .11\textheight,
    scale only axis,
    xmin = 100,
    xmax = 1000000,
    ymin = 2e-5,
    ymax = 1e+2,
    xminorticks = false,
    yminorticks = false,
    xlabel = {\normalsize Frequency $\omega$ (rad/s)},
    ylabel = {\normalsize Magnitude},
    ylabel style   = {yshift = -.3em},
    scaled x ticks = false]
    
    \addplot[sovf1] table[x index = 0, y index = 4] {\tableERR};
    \addplot[sovf2] table[x index = 0, y index = 4] {\tableERRV};
  \end{semilogyaxis}
\end{tikzpicture}%
  \tikzexternaldisable%

    \caption{Pointwise relative errors.}
    \label{fig:sovf_btterfly_relerr}
  \end{subfigure}
  \vspace{0\baselineskip}

  \tikzexternalenable%
  \tikzsetnextfilename{sovf_butterfly_legend}%
  \begin{tikzpicture}  
  \begin{axis}[%
    hide axis,
    xmin = 0,
    xmax = 1,
    ymin = 0,
    ymax = 1,
    scale only axis,
    width = 1mm,
    legend columns = -1, 
    legend style = {
      at     = {(0,0)},
      anchor = center,
      /tikz/every even column/.append style = {column sep = 0.5cm}}]
    
    \addplot[fom] coordinates {(0,0)};
    \addlegendentry{Original data};
    
    \addplot[sovf1] coordinates {(0,0)};
    \addlegendentry{\textsf{SOVF1}};
    
    \addplot[sovf2] coordinates {(0,0)};
    \addlegendentry{\textsf{SOVF2}};
  \end{axis}
\end{tikzpicture}%
  \tikzexternaldisable%

  \caption{Results for the butterfly gyroscope data.}
  \label{fig:sovf_butterfly}
\end{figure}

\subsection{Butterfly gyroscope example}
\label{sec:ButtGyro}

First, we present the results for the butterfly gyroscope model as shown in
\Cref{fig:sovf_butterfly}.
For the given data, we have used the two proposed approaches to learn
structure-preserving models of order $r = 8$.
While \textsf{SOVF2} converges up to numerical accuracy, this is not the case
for \textsf{SOVF1}.
However, the denominator in \textsf{SOVF1} converges reasonably close to
the value one.
Consequently, we have simply considered the mechanical system associated with
the rational function in the numerator of~\cref{eqn:partial_bary_Hr} and
ignored the denominator entry altogether.

As it can be easily observed in \Cref{fig:sovf_btterfly_tf},
\textsf{SOVF1} accurately approximates the given data over the full frequency
range.
A minor exception is given by the right limit of the frequency interval, where
the transfer function of \textsf{SOVF1} slightly deviates from the data.
On the other hand, \textsf{SOVF2} lacks this good approximation behavior as it
is illustrated in \Cref{fig:sovf_btterfly_relerr}.
The accuracy of this approach is at least two orders of magnitude worse than
that of \textsf{SOVF1}.
Still, \textsf{SOVF2} yields a reasonable approximation for the low frequency
range.
We have observed that \textsf{SOVF2} typically introduces poles close
to the imaginary axis, and has better approximation quality in the large
magnitude range.


\begin{figure}[t]
  \centering
  \begin{subfigure}[b]{\linewidth}
    \centering
  \tikzexternalenable%
  \tikzsetnextfilename{sovf_fishtail_tf}%
  \begin{tikzpicture}
  \pgfplotstableread{graphics/data/sovf1_fishtail.dat}\tableTF
  \pgfplotstableread{graphics/data/sovf2_real_fishtail.dat}\tableTFV
  
  \begin{semilogyaxis}[%
    width  = .75\linewidth,
    height = .11\textheight,
    scale only axis,
    xmin = 0,
    xmax = 1000,
    ymin = 1e-8,
    ymax = 3e-3,
    xminorticks = false,
    yminorticks = false,
    xlabel = {\normalsize Frequency $\omega$ (rad/s)},
    ylabel = {\normalsize Magnitude},
    ylabel style   = {yshift = -.3em},
    scaled x ticks = false]
    
    \addplot[fom] table[x index = 0, y index = 1] {\tableTF};
    \addplot[sovf1] table[x index = 0, y index = 2] {\tableTF};
    \addplot[sovf2] table[x index = 0, y index = 2] {\tableTFV};
  \end{semilogyaxis}
\end{tikzpicture}%
  \tikzexternaldisable%

    \caption{Transfer functions.}
  \end{subfigure}
  \vspace{0\baselineskip}
  
  \begin{subfigure}[b]{\linewidth}
    \centering
  \tikzexternalenable%
  \tikzsetnextfilename{sovf_fishtail_relerr}%
  \begin{tikzpicture}
  \pgfplotstableread{graphics/data/sovf1_fishtail.dat}\tableERR
  \pgfplotstableread{graphics/data/sovf2_real_fishtail.dat}\tableERRV
  
  \begin{semilogyaxis}[%
    width  = .75\linewidth,
    height = .11\textheight,
    scale only axis,
    xmin = 0,
    xmax = 1000,
    ymin = 5e-5,
    ymax = 5e-0,
    xminorticks = false,
    yminorticks = false,
    xlabel = {Frequency $\omega$ (rad/s)},
    ylabel = {Magnitude},
    ylabel style   = {yshift = -.3em},
    scaled x ticks = false]
    
    \addplot[sovf1] table[x index = 0, y index = 4] {\tableERR};
    \addplot[sovf2] table[x index = 0, y index = 4] {\tableERRV};
  \end{semilogyaxis}
\end{tikzpicture}%
  \tikzexternaldisable%

    \caption{Pointwise relative errors.}
    \label{fig:sovf1_fishtail_relerr}
  \end{subfigure}
  \vspace{0\baselineskip}

  \tikzexternalenable%
  \tikzsetnextfilename{sovf_fishtail_legend}%
  \begin{tikzpicture}  
  \begin{axis}[%
    hide axis,
    xmin = 0,
    xmax = 1,
    ymin = 0,
    ymax = 1,
    scale only axis,
    width = 1mm,
    legend columns = -1, 
    legend style = {
      at     = {(0,0)},
      anchor = center,
      /tikz/every even column/.append style = {column sep = 0.5cm}}]
    
    \addplot[fom] coordinates {(0,0)};
    \addlegendentry{Original data};
    
    \addplot[sovf1] coordinates {(0,0)};
    \addlegendentry{\textsf{SOVF1}};
    
    \addplot[sovf2] coordinates {(0,0)};
    \addlegendentry{\textsf{SOVF2}};
  \end{axis}
\end{tikzpicture}%
  \tikzexternaldisable%

  \caption{Results for the artificial fishtail data.}
  \label{fig:sovf1_fishtail}
\end{figure}

\subsection{Artificial fishtail example}
\label{sec:ArtFish}

We now present results for the artificial fishtail example, as shown in
\Cref{fig:sovf1_fishtail}.
We have used both approaches and constructed structure-preserving learned
models of order $r = 10$. As in the first example presented in
\Cref{sec:ButtGyro},
\textsf{SOVF1} provides a high-fidelity approximation over the full frequency
interval and outperforms \textsf{SOVF2}.
As shown in \Cref{fig:sovf1_fishtail_relerr}, the \textsf{SOVF2}
approximation has a large mismatch around $580$\,rad/s, where the relative
approximation error is the highest.
Like in the previous example, \textsf{SOVF2} to provides an accurate
approximation of the main dominant peak, i.e., the one located around
$100$\,rad/s.
Further details on the numerical results and convergence of the applied
methods can be found in the accompanying code package.


\section{Conclusions}%
\label{sec:conclusions}

We have proposed two new approaches for data-driven modeling of modally damped
mechanical systems by developing  structure-preserving vector fitting
formulations.
We have revised the barycentric formula to represent structured transfer
functions, and have shown that the structure of the original model is
automatically preserved in the reduced one.
The two approaches have been applied to two benchmark models and the
preliminary results are promising.
The method corresponding to the partially structured transfer function
formulation has been proven especially accurate and reliable in both test
cases.
A more thorough investigation is needed to explain the accuracy miss-matches
encountered in the \textsf{SOVF2} formulation.
Extending the analysis and numerical algorithms to MIMO problems and employing
the proposed structured barycentric forms to develop a AAA-like framework are
natural next steps.


\section*{Acknowledgments}%
\addcontentsline{toc}{section}{Acknowledgments}

Gosea and Werner, while he was at Max Planck Institute Magdeburg, have been
supported in parts by the German Research Foundation (DFG) Research Training
Group 2297 \textquotedblleft{}MathCoRe\textquotedblright{}, Magdeburg.
Gugercin was supported in parts by the National Science Foundation under
Grant No. DMS-1923221.


\addcontentsline{toc}{section}{References}
\bibliographystyle{plainurl}
\bibliography{bibtex/myref}

\end{document}